\documentclass[reqno,dvips,12pt]{article} 
\usepackage{epsfig,graphicx,amsmath,amssymb,amsfonts}
cm \textwidth=16.0 true cm \emergencystretch=10pt

\newtheorem{theorem}{\qquad Theorem}
\newtheorem{lemma}{\qquad Lemma} 
\newcommand{\abs}[1]{\left\vert{#1}\right\vert}
\newcommand{\e}{\varepsilon}

\DeclareMathOperator{\re}{re}

\newcommand{\cA}{{\cal A}}
\newcommand{\cM}{{\cal M}}

\newcommand{\cN}{{\cal N}}
\newcommand{\cP}{{\cal P}}

\newcommand{\const}{\mathop{\rm const}\nolimits}

\newcommand{\nn}{\nonumber}

\newcommand{\pa}{\partial}

\begin{document}

\title {\textbf{Negative Dimension  in General and Asymptotic Topology}}

\author{\textbf{V.P.Maslov}\thanks{Moscow State University,
Physics Department, v.p.maslov@mail.ru}}
\date{ }

\maketitle
\begin{abstract}
We introduce  the notion of negative  topological dimension  and
the notion of weight for the asymptotic topological dimension.
Quantizing of spaces of negative dimension  is applied to
linguistic statistics.
\end{abstract}

\textit{Recently, Yu.~I.~Manin has presented his considerations
for the density of spaces of negative dimension}~\cite{Manin}.

\textbf{ 1.}  Let us consider the simplest examples of (Haar)
measures in the general case for the $n$-dimensional space. Let
$S_n$ be an $n$-dimensional ball of radius~$r$. In the spherical
coordinates, the volume $\mu(S_n)$ of the ball is equal to $const
\int_0^1 r^{n-1} dr= const \  r^n$. Here $r^{n-1}$ stands for the
density.

In the sense of the Fourier transform, the multiplication by a
coordinate is dual to the corresponding derivation. Therefore, we
can speak of dual $n$ times differentiable functions in the
Sobolev space~$W_2^n$.  Dirac distinguished between the left and
right components ''bra'' and ''ket'' in the ''bracket'' inner
product. The ``dual'' space of this space according to Dirac is
the space~$W_2^{-n}$ of Sobolev distributions (generalized
functions).

In the same way we can define the functions in $W_2^s$ by the
``inner product,'' where $s$ is a positive noninteger number, and
the space $W_2^{-s}$ as the ``inner product'' conjugate
to~$W_2^s$.

One can similarly proceed with the density (or the weight) $r^s$
and $r^{-s}$, by using, for instance, the Riesz kernel or the
Bessel potential to represent functions in~$W_2^s$.

Let us present an example of a space (of noninteger positive
dimension) equipped with the Haar measure $r^\sigma$, where $0
\leq\sigma \leq 1$.

On the closed interval $0\leq x\leq 1$ there is a scale $0\leq
\sigma\leq 1$ of Cantor dust with the Haar measure equal to
$x^{\sigma}$ for any interval $(0,x)$ similar to the entire given
set of the Cantor dust. The direct product of this scale by the
Euclidean cube of dimension~$k-1$ gives the entire scale
$k+\sigma$, where $k \in \Bbb Z$ and~$\sigma \in (0,1)$.

\textit{General definition of spaces of negative dimension.} Let
$M_{t_0}$ be a compactum, of Hausdorff dimension~$t_0$, which is
an element of a $t$-parameter scale of mutually embedded
compacta, $0< t < \infty$. Two scales of this kind are said to be
{\em equivalent\/} with respect to the compactum $M_{t_0}$ if all
compacta in these scales coincide for any $t\geq t_0$. We say
that the compactum $M_{t_0}$ is a {\em hole\/} in this equivalent
set of scales and the number $-t_0$ is the {\em negative
dimension\/} of this equivalence class.

We consider the space of negative dimension $-D=-k-\sigma$ with
respect to the given above scale.

\textbf{2.} As in \cite{NegDimen}, the values of the random
variable $x_1, \dots, x_s$ are ordered in absolute value. Some of
the numbers $x_1,\dots,x_s$ may coincide. Then these numbers are
combined adding the corresponding ''probabilities'', i.e., the
ratio of the number of ``hits'' at~$x_i$ to the general number of
trials. The number of equal $x_i: x_i=x_{i+1}= \dots = x_{i+k}$
will be called the multiplicity $q_i$ of the value $x_i$. In our
consideration, both the number of trials $N$ and~$s$ tend to
infinity.

Let $N_i$ be the number of ''appearances'' of the value $x_i: \
x_i < x_{i+1}$, then
\begin{equation}
\sum^s_{i=1} \frac{N_i}{N} x_i=M, \label{Zipf1}
\end{equation}
where $M$ is the mathematical expectation.

The cumulative probability $\cP_k$ is the sum of the first~$k$
probabilities in the sequence $x_i$: $\cP_k=\frac 1N \sum_{i=1}^k
N_i$, where $k<s$. We denote $NP_k=B_k$.

If all the variants for which
\begin{equation}\label{A}
\sum_{i=1}^s N_i = N
\end{equation}
and
\begin{equation}\label{B}
\sum_{i=1}^s N_i x_i \leq E, \ \ E=MN\leq N \overline{x},
\end{equation}
where $\overline{x}=\frac{\sum_{i=1}^s q_i x_i}{Q}$,
$Q=\sum_{i=1}^s q_i$, are equivalent (equiprobable), then
\cite{MatZamTheor,NelinSred,NoPredp} the majority of the variants
will accumulate near the following dependence of the ''cumulative
probability'' $B_l\{N_i\}=\sum_{i=1}^l N_i$,
\begin{equation}
\sum_{i=1}^l N_i= \sum_{i=1}^l \frac{q_i}{e^{\beta'x_i-\nu'}-1},
\label{Zipf2}
\end{equation}
where $\beta'$ and $\nu'$ are determined by the conditions
\begin{equation}\label{Zipf2a}
B_s=N,
\end{equation}
\begin{equation}\label{Zipf2a'}
\sum_{i=1}^s \frac{q_i x_i}{e^{\beta' x_i-\nu'}-1}=E,
\end{equation}
as $N \to \infty$ and $s \to \infty$.

We introduce the notation: $\cM$ is the set of all sets $\{N_i\}$
satisfying conditions~(\ref{A}) and~(\ref{B}); \ $\cN\{\cM\}$ is
the number of elements of the set~$\cM$.

\begin{theorem} \label{theor1}
Suppose that all the variants of sets $\{N_i\}$ satisfying the
conditions ~(\ref{A}) and ~(\ref{B}) are equiprobable. Then the
number of variants $\cN$ of sets $\{N_i\}$ satisfying
conditions~(\ref{A}) and~(\ref{B}) and the additional relation
\begin{equation} \label{theorema1}
|\sum^l_{i=1} N_i - \sum^l_1\frac{q_i}{e^{\beta'
x_i-\nu'}-1}|\geq  N^{(3/4+\varepsilon)}
\end{equation}
is less than $\frac{c_1 \cN\{\cM\}}{N^m}$ (where~$c_1$ and~$m$
are any arbitrary numbers, $\sum_{i=1}^l q_i \geq\varepsilon Q$,
and $\varepsilon$ is arbitrarily small).
\end{theorem}

{\bf{\qquad Proof of Theorem 1.}}

Let $\cA$ be a subset of $\cM$ satisfying the condition
$$
|\sum_{i=l+1}^s N_i - \sum_{i=l+1}^s \frac{q_i} {e^{\beta
x_i-\nu}-1}|\leq \Delta;
$$
$$
|\sum_{i=1}^l N_i-\sum_{i=1}^l \frac{q_i}
{e^{\beta'x_i-\nu'}-1}|\leq \Delta,
$$
where $\Delta$, $\beta$, $\nu$ are some real numbers independent
of~$l$.

We denote
$$
|\sum_{i=l+1}^s N_i-\sum_{i=l+1}^s \frac{q_i} {e^{\beta
x_i-\nu}-1}| =S_{s-l};
$$
$$
|\sum_{i=1}^l N_i-\sum_{i=1}^l \frac{q_i} {e^{\beta'x_i-\nu'}-1}|
=S_l.
$$

Obviously, if $\{N_i\}$ is the set of all sets of integers on the
whole, then
\begin{equation}\label{Proof1}
\cN\{\cM \setminus \cA\} = \sum_{\{N_i\}} \Bigl(
\Theta(E-\sum_{i=1}^s N_ix_i) \delta_{(\sum_{i=1}^s N_i),N}
\Theta(S_l-\Delta)\Theta (S_{s-l}-\Delta)\Bigr),
\end{equation}
where $\sum N_i=N$.

Here the sum is taken over all integers $N_i$, $\Theta(x)$ is the
Heaviside function, and $\delta_{k_1,k_2}$ is the Kronecker
symbol.

We use the integral representations
\begin{eqnarray}
&&\delta_{NN'}=\frac{e^{-\nu N}}{2\pi}\int_{-\pi}^\pi d\varphi
e^{-iN\varphi} e^{\nu N'}e^{i N'\varphi},\label{D7}\\
&&\Theta(y)=\frac1{2\pi i}\int_{-\infty}^\infty
dx\frac1{x-i}e^{\beta y(1+ix)}.\label{D8}
\end{eqnarray}

Now we perform the standard regularization. We replace the first
Heaviside function~$\Theta$ in~(\ref{Proof1}) by the continuous
function
\begin{equation}
\Theta_{\alpha}(y) =\left\{
\begin{array}{ccc}
0 &\mbox{for}& \alpha > 1, \  y<0 \nn \\
1-e^{\beta y(1-\alpha)} &\mbox{for}& \alpha > 1,\  y \geq 0,
\label{Naz1}
\end{array}\right.
\end{equation}
\begin{equation}
\Theta_{\alpha}(y) =\left\{
\begin{array}{ccc}
e^{\beta y(1-\alpha)} &\mbox{for}&\alpha < 0, \ y<0 \nn \\
1 &\mbox{for}& \alpha < 0, \ y \geq 0, \label{Naz2}
\end{array}\right.
\end{equation}
where $\alpha \in (-\infty,0) \cup (1, \infty)$ is a parameter,
and obtain
\begin{equation}\label{proof2}
\Theta_\alpha(y) = \frac1{2\pi i} \int_{-\infty}^{\infty}
e^{\beta y(1+ix)} (\frac 1{x-i} - \frac 1{x-\alpha i}) dx.
\end{equation}

If  $\alpha > 1$, then $\Theta(y)\leq \Theta_{\alpha}(y)$.

Let $\nu <0$. We substitute~(\ref{D7}) and~(\ref{D8})
into~(\ref{Proof1}),  interchange the integration and summation,
then pass to the limit as $\alpha \to \infty$ and obtain the
estimate
\begin{eqnarray}
&&\cN\{\cM \setminus   \cA\} \leq \nn \\
&&\leq \Bigl|\frac{e^{-\nu N+\beta E }}{i(2\pi)^2}\int_{-\pi}^\pi
\bigl[ \exp(-iN\varphi)
\sum_{\{N_j\}}\exp\bigl\{-\beta\sum_{j=1}^s N_jx_j+(i\varphi+\nu)
\sum_{j=1}^s N_j\bigr\}\bigr]\ d\varphi \times \nn \\
&& \times \Theta(S_l -\Delta)\Theta(S_{s-l}-\Delta)\Bigr|,
\end{eqnarray}
where $\beta$ and $\nu$ are real parameters such that the series
converges for them.

To estimate the expression in the right-hand side, we bring the
absolute value sign inside the integral sign and then inside the
sum sign, integrate over $\varphi$, and obtain
\begin{eqnarray}
&&\cN\{\cM \setminus \cA\} \leq \frac{e^{-\nu N+\beta E }}{2\pi}
\sum_{\{N_i\}}\exp\{-\beta\sum_{i=1}^sN_ix_i+\nu
\sum_{i=1}^s N_i\}\times \nn \\
&& \times\Theta (S_l-\Delta)\Theta (S_{s-l}-\Delta).
\end{eqnarray}

We denote
\begin{equation}\label{D9a}
Z(\beta,N)=\sum_{\{N_i\}} e^{-\beta\sum_{i=1}^s N_ix_i},
\end{equation}
where the sum is taken over all $N_i$ such that $\sum_{i=1}^s
N_i=N$,
$$
\zeta_l(\nu,\beta)= \prod_{i=1}^{l} \xi_i\left(\nu,\beta\right);
\zeta_{s-l}(\nu,\beta)= \prod_{i=l+1}^{s}
\xi_i\left(\nu,\beta\right);
$$
$$
\quad \xi_i(\nu,\beta)= \frac{1}{(1-e^{\nu-\beta x_i})^{q_i}},
\qquad i=1,\dots,l.
$$

It follows from the inequality for the hyperbolic cosine
$\cosh(x)=(e^x+e^{-x})/2$ for $|x_1| \geq \delta; |x_2| \geq
\delta$:
\begin{equation}
\cosh(x_1)\cosh(x_2)= \cosh(x_1+x_2) + \cosh(x_1 - x_2) >
\frac{e^\delta}{2} \label{D33}
\end{equation}
that the inequality
\begin{equation}
\Theta(S_{s-l}-\Delta) \Theta(S_{l}-\Delta)\le e^{-c\Delta}
\cosh\Bigl(c\sum_{i=1}^{l} N_i
-c\phi_l\Bigr)\cosh\Bigl(c\sum_{i=l+1}^{s} N_i
-c\overline{\phi}_{s-l}\Bigr), \label{D34}
\end{equation}
where
$$
\phi_l= \sum_{i=1}^l \frac{q_i}{e^{\beta'x_i-\nu'}-1}; \qquad
\overline{\phi}_{s-l}= \sum_{i=l+1}^s \frac{q_i}{e^{\beta
x_i-\nu}-1},
$$
holds for all positive $c$ and~$\Delta$.

We obtain
\begin{eqnarray}
&&\cN\{\cM \setminus \cA\} \leq  e^{-c\Delta} \exp\left(\beta
E-\nu N\right) \times \nn \\
&& \times \sum_{\{N_i\}}\exp\{-\beta\sum_{i=1}^l
N_ix_i+\nu\sum_{i=1}^l N_i\} \cosh\left(\sum_{i=1}^{l} c N_i -
c\phi\right) \times \nn \\
&& \times \exp\{-\beta\sum_{i=l+1}^s N_ix_i +\nu \sum_{i=l+1}^s
N_i\} \cosh\Bigl(\sum_{i=l+1}^s c N_i -c\overline{\phi}\Bigr) = \nn \\
&& =e^{\beta E} e^{-c\Delta} \times \nn \\
&&\times \left( \zeta_l(\nu-c,\beta) \exp(-c\phi_{l})
+\zeta_l(\nu+c,\beta)\exp(c\phi_{l})\right) \times \nn \\
&&\times\left(\zeta_{s-l}(\nu-c,\beta)
\exp(-c\overline{\phi}_{s-l})
+\zeta_{s-l}(\nu+c,\beta)\exp(c\overline{\phi}_{s-l})\right).
\label{5th}
\end{eqnarray}

Now we use the relations
\begin{equation}\label{5tha}
\frac {\pa}{\pa\nu}\ln \zeta_l|_{\beta=\beta',\nu=\nu'}\equiv
\phi_l; \frac {\pa}{\pa\nu}\ln
\zeta_{s-l}|_{\beta=\beta',\nu=\nu'}\equiv \overline{\phi}_{s-l}
\end{equation}
and the expansion $\zeta_l(\nu\pm c,\beta)$ by the Taylor formula.
There exists a $ \gamma <1$ such that
$$
\ln(\zeta_l(\nu\pm c,\beta)) =\ln\zeta_l(\nu,\beta)\pm
c(\ln\zeta_l)'_\nu(\nu,\beta)+\frac{c^2}{2}(\ln\zeta_l)^{''}_\nu
(\nu\pm\gamma c,\beta).
$$
We substitute this expansion, use formula~(\ref{5tha}), and see
that $\phi_{\nu,\beta}$ is cancelled.

Another representation of the Taylor formula implies
\begin{eqnarray}
&&\ln\left(\zeta_l(\nu+c,\beta)\right)=
\ln\left(\zeta_l(\beta,\nu)\right)+
\frac{c}\beta\frac{\pa}{\pa\nu}\ln\left(\zeta_l(\beta,\nu)\right)+\nn\\
&&+\int_{\nu}^{\nu+c/\beta}d\nu' (\nu+c/\beta-\nu')
\frac{\pa^2}{\pa\nu'^2}\ln\left(\zeta_l(\beta,\nu')\right).\label{CC1}
\end{eqnarray}
A similar expression holds for $\zeta_{s-l}$.

From the explicit form of the function $\zeta_l(\beta,\nu)$, we
obtain
\begin{equation}
\frac{\pa^2}{\pa\nu^2}\ln\left(\zeta_l(\beta,\nu)\right)=
\beta^2\sum_{i=1}^{l}
\frac{g_i\exp(-\beta(x_i+\nu))}{(\exp(-\beta(x_i+\nu))-1)^2}\leq
\beta^2Qd, \label{CC2}
\end{equation}
where $d$ is given by the formula
$$
d=\frac{\exp(-\beta(x_1+\nu))}{(\exp(-\beta(x_1+\nu))-1)^2}..
$$
The same estimate holds for $\zeta_{s-l}$.

Taking into account the fact that $\zeta_l\zeta_{s-l}=\zeta_s$,
we obtain the following estimate for $\beta=\beta'$ and
$\nu=\nu'$:
\begin{equation}\label{eval1}
\cN\{\cM \setminus \cA\}
\leq\zeta_s(\beta',\nu')\exp(-c\Delta+\frac{c^2}{2}\beta^2Qd)
\exp(E\beta'-\nu'N).
\end{equation}

Now we express $\zeta_s(\nu',\beta')$ in terms $Z(\beta,N)$. To
do this, we prove the following lemma.

\begin{lemma}
Under the above assumptions, the asymptotics of the integral
\begin{equation}\label{lemma_1}
Z(\beta,N) = \frac{e^{-\nu N}}{2\pi}\int_{-\pi}^\pi
 d\alpha e^{-iN\alpha}\zeta_s(\beta,\nu+i\alpha)
\end{equation}
has the form
\begin{equation}\label{lemma_2}
Z(\beta,N) = C e^{-\nu N} \frac{\zeta_s(\beta,\nu)}{|(\partial^2
\ln\zeta_s(\beta,\nu))/ (\partial^2\nu)|} (1+O(\frac 1N)),
\end{equation}
where $C$ is a constant.
\end{lemma}

We have
\begin{equation}
 Z(\beta,N) = \frac{e^{-\nu N}}{2\pi}\int_{-\pi}^\pi
 e^{-iN\alpha}\zeta_s(\beta,\nu+i\alpha)\,d\alpha
 =\frac{e^{-\nu N}}{2\pi}\int_{-\pi}^\pi e^{NS(\alpha,N)} d\alpha ,\label{D15}
\end{equation}
where
\begin{equation}\label{qq}
    S(\alpha,N) = -i\alpha+ \ln \zeta_s (\beta, \nu +i\alpha)
    = -i\alpha - \sum_{i=1}^s q_i\ln [1-e^{\nu+i\alpha-\beta
    x_i}].
\end{equation}
Here $S$ depends on $N$, because $s$, $x_i$, and $\nu$ also
depend on~$N$; the latter is chosen so that the point $\alpha=0$
be a stationary point of the phase~$S$, i.e., from the condition
\begin{equation}\label{qq1}
 N=\sum_{i=1}^s\frac{q_i}{e^{\beta
 x_i-\nu}-1}.
\end{equation}
We assume that $a_1N \leq s \leq a_2N$, $a_1,
a_2=\operatorname{const}$, and, in addition, $0\le x_i\le B$ and
$B=\operatorname{const}$, $i=1,\dots,s$.
 If these conditions are satisfied
in some interval $\beta\in[0,\beta_0]$ of the values of the
inverse temperature, then all the derivatives of the phase are
bounded, the stationary point is nondegenerate, and the real part
of the phase outside a neighborhood of zero is strictly less than
its value at zero minus some positive number. Therefore,
calculating the asymptotics of the integral, we can replace the
interval of integration $[-\pi,\pi]$ by the interval $[-\e,\e]$.
In this integral, we perform the change of variable
\begin{equation}\label{qqq}
  z=\sqrt {S(0,N)-S(\alpha,N)}.
\end{equation}
This function is holomorphic in the disk $\abs{\alpha}\le\e$ in
the complex $\alpha$-plane and has a holomorphic inverse for a
sufficiently small~$\e$. As a result, we obtain
\begin{equation}\label{qqq1}
    \int_{-\e}^\e e^{NS(\alpha,N)}
    d\alpha=e^{NS(0,N)}\int_\gamma e^{-Nz^2}f(z)\,dz,
\end{equation}
where the path $\gamma$ in the complex $z$-plane is obtained from
the interval $[-\e,\e]$ by the change~\eqref{qqq} and
\begin{equation}\label{qqq3}
    f(z)=\left(\frac{\partial\sqrt {S(0,N)-S(\alpha,N)}}
    {\partial\alpha}\right)^{-1}\bigg|_{\alpha=\alpha(z)}.
\end{equation}
For a small~$\e$ the path $\gamma$ lies completely inside the
double sector $\re(z^2)>c(\re z)^2$ for some $c>0$; hence it can
be ``shifted'' to the real axis so that the integral does not
change up to terms that are exponentially small in~$N$. Thus,
with the above accuracy, we have
\begin{equation}\label{qqq4}
  Z(\beta,N) =  \frac{e^{-\nu N}}{2\pi}\int_{-\e}^\e e^{-Nz^2}f(z)\,dz.
\end{equation}
Since the variable $z$ is now real, we can assume that the
function $f(z)$ is finite (changing it outside the interval of
integration), extend the integral to the entire axis (which again
gives an exponentially small error), and then calculate the
asymptotic expansion of the integral expanding the integrand in
the Taylor series in~$z$ with a remainder. This justifies that
the saddle-point method can be applied to the above integral in
our case.

\begin{lemma}
The quantity
\begin{equation}\label{qqq5}
\frac{1}{\cN(\cM)} \sum_{\{N_i\}} e^{-\beta\sum_{i=1}^s N_ix_i},
\end{equation}
where $\sum N_i =N$ and $x_iN_i\leq E-N^{1/2+\varepsilon}$, tends
to zero faster than $N^{-k}$ for any $k$, $\varepsilon>0$.
\end{lemma}

We consider the point of minimum in $\beta$ of the right-hand
side of ~(\ref{5th}) with $\nu(\beta,N)$ satisfying the condition
$$
\sum \frac{q_i}{e^{\beta x_i-\nu(\beta,N)}-1} =N.
$$
It is easy to see that it satisfies condition~(\ref{Zipf2a}). Now
we assume that the assumption of the lemma is not satisfied.

Then for $\sum N_i=N$,  $\sum x_i N_i\geq E-N^{1/2+\varepsilon}$,
we have
$$
e^{\beta E}\sum_{\{N_i\}} e^{-\beta\sum_{i=1}^s N_ix_i}\geq
e^{(N^{1/2}+\varepsilon)\beta}.
$$
Obviously, $\beta\ll \frac{1}{\sqrt{N}}$ provides a minimum
of~(\ref{5th}) if the assumptions of Lemma~1 are satisfied, which
contradicts the assumption that the minimum in~$\beta$ of the
right-hand side of~(\ref{5th}) is equal to~$\beta'$.

We set $c=\frac\Delta{N^{1+\alpha}}$ in formula~(\ref{eval1})
after the substitution~(\ref{lemma_2}); then it is easy to see
that the ratio
$$
\frac{\cN(\cM \setminus\cA)}{\cN(\cM)}\approx \frac 1{N^m},
$$
where $m$ is an arbitrary integer, holds for
$\Delta=N^{3/4+\varepsilon}$. The proof of the theorem is
complete.

We will prove a cumulative formula in which the densities coincide
in shape with the Bose--Einstein distribution. The difference
consists only in that, instead of the set $\lambda_n$ of random
variables or eigenvalues of the Hamiltonian operator, the
Bose--Einstein formula contains some of their averages over the
cells \cite{Landau1}. In view of the theorem given below,  one can
proof that  the $\varepsilon_i$, which are averages of the energy
$\lambda_k$ at the $i$th cell, are nonlinear averages in the
sense of Kolmogorov~\cite{NelinSred}.

\textbf{3.} Now we consider the notion of the lattice dimension.

We consider a straight line, a plane, and a three-dimensional
space. We separate points $i=0,1,2, \dots$ on the line and points
$x=i=0,1,2,\dots$, \ $y=j=0,1,2, \dots$ on the coordinate axes
$x,y$ on the plane. We associate this set of points $(i,j)$ with
the points on the straight line (with the positive integers
$l=1,2 \dots$) up to the quantum constant~$\chi$ of the lattice.

According to M.~Gromov's definition~\cite{Gromov}, the asymptotic
(topological) dimension of this lattice is equal to two.

We associate each point with a pair of points~$i$ and~$j$
according to the rule $i+j=l$. The number of such points $n_l$ is
equal to $l+1$. In addition, we assume that $z=k=0,1,2,\dots$ on
the axis, i.e., we set $i+j+k=l$. In this case, the number of
points $q_l$ is equal to
$$
q_l =\frac{(l+1)(l+2)}{2}.
$$

If we set $\lambda_i=l$ in formula~(\ref{Zipf2}), then, in the
three-dimensional case, each~$i$ is associated with
$\frac{(l+1)(l+2)}{2}$ of mutually equal $x_l=l$ (these are the
multiplicities or the $q_l-$hold degeneracies of the spectrum of
the oscillator). Formula~(\ref{Zipf2}) in this special case
becomes
\begin{equation}\label{3.4}
N_l= const \sum_{i=0}^l \frac{(i+1)(i+2)}{2(e^{\beta i-\nu}-1)};
\end{equation}
\begin{equation}\label{3.5}
\Delta N_i= const \frac{(i+1)(i+2)}{2(e^{\beta
i-\nu}-1)}\Delta_i, \ \ \Delta_i=1,
\end{equation}
\begin{equation}\label{3.5a}
\Delta E_i=\const \frac{i(i+1)(i+2)}{2(e^{\beta i-\nu}-1)}\Delta_i
\end{equation}
for large~$i$, $\frac{\Delta_i}{i} \to 0$,
\begin{equation}\label{3.6}
dE= const \frac{\omega^3 d \omega}{e^{\beta\omega}-1}; \ \beta=
\frac hT
\end{equation}
(cf. formula~(60.4) in~\cite{Landau}).

Thus, we obtain a somewhat sharper version of the famous Planck
formula for the radiation of a black body.

For the $D$-dimensional case, it is easy to verify that the
sequence of weights (multiplicities) of the number of versions
$i= \sum_{k=1}^D m_k$, where $m_k$ are arbitrary positive
integers, has the form of the binomial coefficient
\begin{equation}\label{3.7}
q_i(D) = const\frac{(i+D-1)!}{i!D!},
\end{equation}
where the constant depends on~$D$.

Thus, for any $D$, formula~(\ref{Zipf2}) has the form
\begin{equation}\label{3.8}
N_l = const \sum_{i=1}^l  \frac{q_i (D)}{e^{\beta i}-1}.
\end{equation}

For the positive integers, we have a sequence of weights~$q_i$
(or, simply, a weight) of the form~(\ref{3.7}).

Our weight series can easily be continued to an arbitrary case by
replacing the factorials with the $\Gamma$-functions; in this
case, we assume that~$D$ is negative.

This is the negative topological dimension (the hole dimension)
of the quantized space (lattice).

If $D>1$, then, as $i\to\infty$, a condensation of a sufficiently
small perturbation occurs in the spectrum of the oscillator and
the multiplicities split, i.e., the spectrum becomes denser
as~$i$ increases. The fact that~$D$ is negative means that there
is strong rarefaction in the spectrum as $i\to\infty$ (the
constant in formula~(\ref{3.8}) must be sufficiently large).

For non-positive integer~$D$, the terms $i=0,1,2,3, \dots, -D$
become infinite. This means that they are very large in the
experiment, which permits determining the lattice negative
dimension corresponding to a given problem. We note that a new
condensate occurs, which is possible for small~$\beta$.

\textbf{4.} Now we consider an example of negative asymptotic
dimension in linguistic statistics.

To each word a frequency dictionary assigns its occurrence
frequency (i.e., the number of occurrences) of this word in the
corresponding corpus of texts. There may be several words with
same frequency.

There is an analogy between the Bose particles at the energy
level of an oscillator $\lambda_i=i$ and the words with the
occurrence frequency (occurrence number)~$i$, namely, the words
with the same occurrence frequency can be ordered in an arbitrary
way, say, alphabetically, inversely alphabetically, or in any
other order. The indexing of the ranks (the indices) of words
within the family of a given occurrence frequency (occurrence
number) is arbitrary. In this sense, the words are
indistinguishable and are distributed according to the Bose
statistic.

However, there is a difference between the approaches under
consideration. In the frequency dictionary one evaluates the
occurrence frequency of every word and then orders the words,
beginning with the most frequently occurring words.

When there were no computers, it was difficult for a person to
evaluate the number of words with equal occurrence frequency. By
looking at a page as if it were a picture, a person can determine
a desired word on this page by its graphical form at every place
of occurrence of the word. In this case, the person looks at a
page of the text as if it were a photo, without going into the
meaning. Similarly, if a person looks for a definite name in a
long list of intrants who had entered a college, this person
finds (or does not find) the desired name by eyes rather than
reads all the names one after another.

An eye gets into the way of recognizing the desired image, and
this ability intensifies as the viewed material increases: the
more pages the eye scans, the less is the difficulty in finding
the desired graphical form. Therefore, under a manual counting,
it was simpler to recognize the desired word on a page without
reading the text and to cross it out by a pencil, simultaneously
counting the number of occurrences of the word. This procedure is
repeated for any subsequent word, already using the text with the
words crossed out (``holes''), which facilitates the search. In
other words, the procedure is in the recognition of the image of
the given word, similar to the recognition, say, of a desired
mushroom in a forest without sorting out all the plants on the
soil one after another. An ordinary computer solves problems of
this kind by exhaustion, whereas a quantum computer
(see~\cite{Belavkin_Mas}) makes this by recognizing the image.

However, for an ordinary computer, the number of operations
needed to find the occurrence frequency of a word is less than
the number of operations needed to find the number of words in
the text with a given occurrence frequency.

One can say that the number of mushrooms we gathered (took away
from the forest) is the number of holes we left in the forest.
Similarly, the words we had ``got out'' from the text in the
above way is an analog of holes rather than particles. Therefore,
the linguists count the rank of words starting from the opposite
end as compared to the starting end which would be used by
physicists. The physicists would count the particles starting
from the lowest level, whereas the holes, the absent electrons,
would be counted from the highest level.

For this reason, the words in a frequency dictionary are
associated with holes rather than particles. Correspondingly, the
dimension in the distribution for  frequency dictionaries  is to
be chosen as a ``hole'' dimension, which is negative.

The number of words encountered only once in the array of texts
is approximately equal to 1/3  of the entire frequency dictionary
which the  number of words equal to $N$. So  as $N\to\infty$ this
is the condensate. It follows from the above that $D = -1$ for the
dictionary. Hence, for $\beta\ll 1$ and $\nu \sim 1$, we have
\begin{equation} \label{rang1}
N_l = const\sum_{i=2}^l \frac{1}{i(i-1)(e^{\beta i-\nu}-1)} \sim
const\int^\omega
\frac{d\omega}{\alpha\omega(\alpha\omega-1)(e^{\beta\alpha\omega-\nu}-1)},
\end{equation}
where $\omega=l$ and $\alpha$ is the scale constant.  If $\omega$
is  finite and  $\beta\ll 1$ the integral may be taken. Therefore
one may set $\beta=0$ for not too large $\omega$.

Bellow we present some graphs of frequency dictionaries obtained
from writings of A. N. Ostrovsky and  N. Berdyaev.

\begin{figure}[h] 
\centering\epsfig{figure=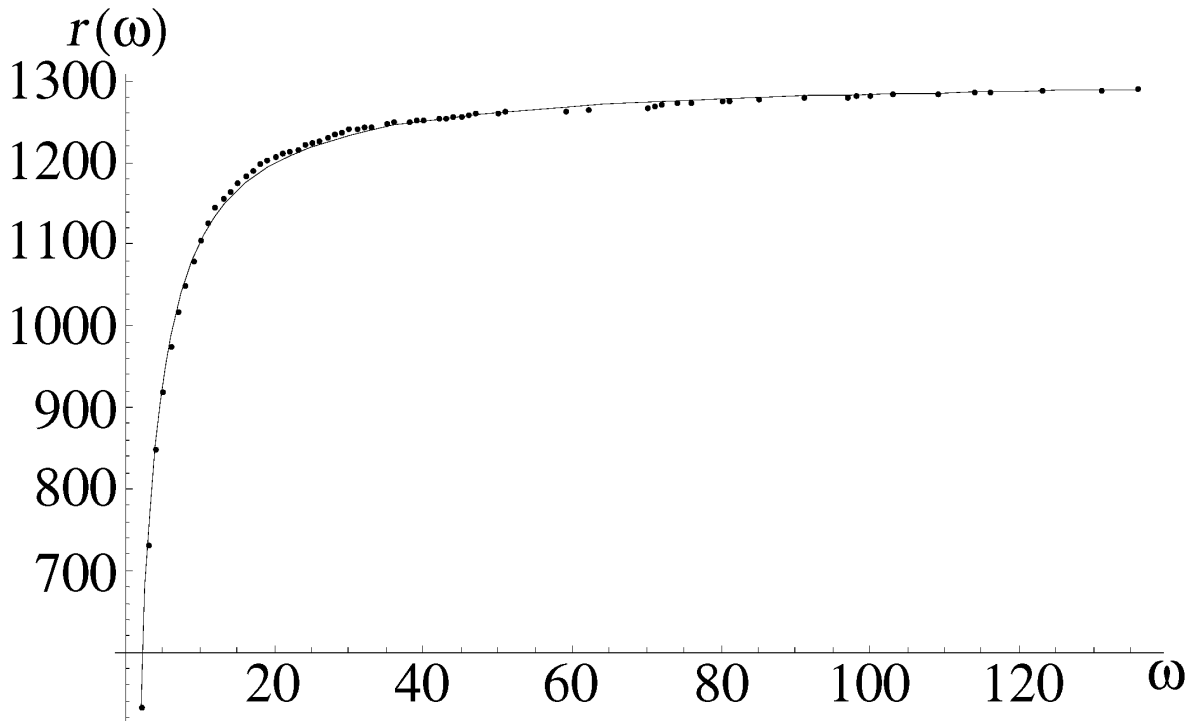,width=\linewidth} \caption{\small
A.N.Ostrovsky "Sin and Trouble Avoid Nobody." [Grekh da beda na
kogo ne zhivet.] Dependence of the inverted rank on the
frequency. Practical and theoretical dependencies. At the bottom,
2 frequencies are cut off: $\omega_{\min} =1$. At the top, 15
frequencies are cut off: $\omega_{\max}= 136$. The approximating
curve is constructed with constants found from the following two
points: $\omega_{1} =2; \omega_{2}= 43, \alpha=0.5.$ Mean
quadratic error: $\sigma = 0.71589$.}
\end{figure}
\begin{figure}[h] 
\centering\epsfig{figure=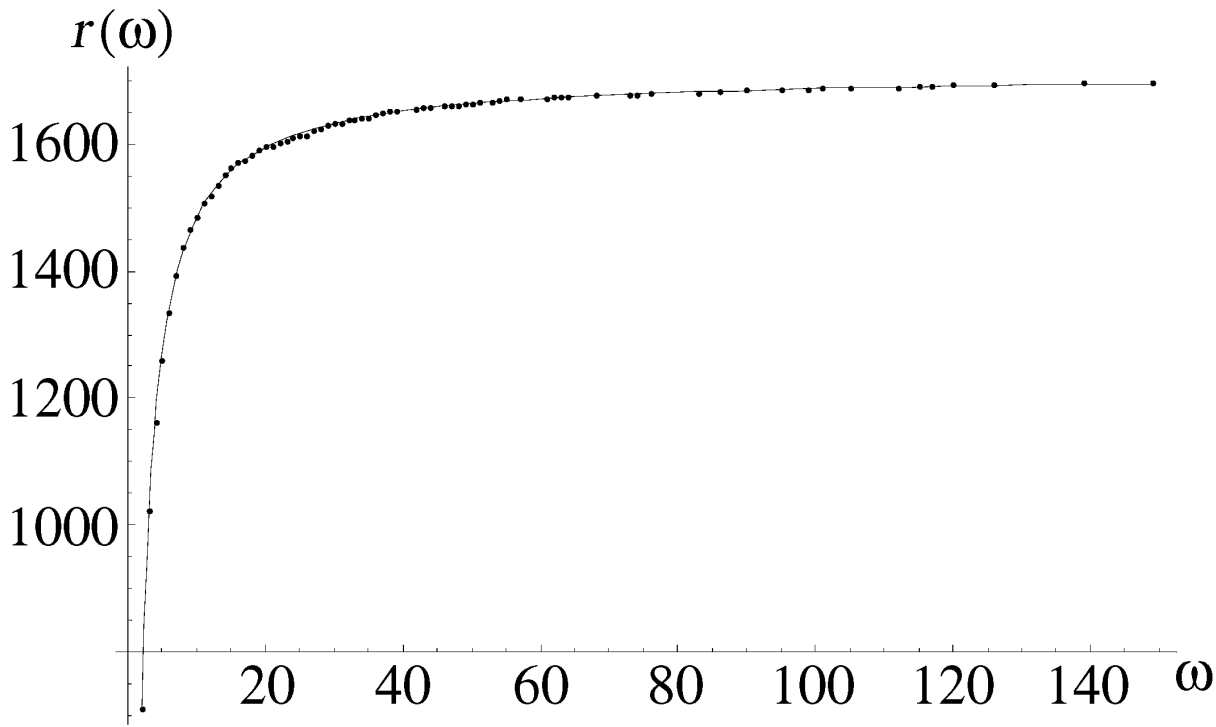,width=\linewidth} \caption{\small
A.N.Ostrovsky "Lucrative Position." [Dokhodnoe mesto.] Dependence
of the inverted rank on the frequency. Practical and theoretical
dependencies. At the bottom, 2 frequencies are cut off:
$\omega_{\min} =1$. At the top, 15 frequencies are cut off:
$\omega_{\max}= 149$. The approximating curve is constructed with
constants found from the following two points: $\omega_{1} =2;
\omega_{2}= 42, \alpha=1.4.$ Mean quadratic error: $\sigma =
0.435321.$}
\end{figure}

\begin{figure}[h] 
\centering\epsfig{figure=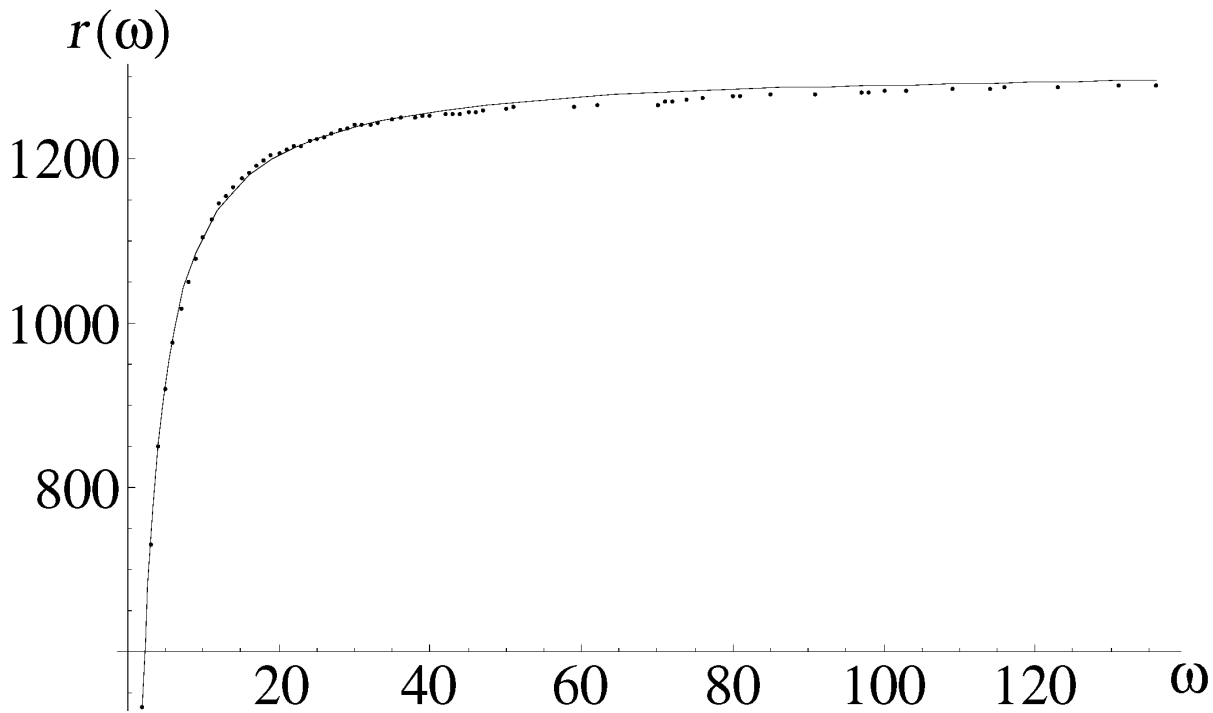,width=\linewidth} \caption{\small
A.N.Ostrovsky "Handsome Gentleman." [Krasavets muzhchina.]
Dependence of the inverted rank on the frequency. Practical and
theoretical dependencies. At the bottom, 2 frequencies are cut
off: $\omega_{\min} =1$. At the top, 15 frequencies are cut off:
$\omega_{\max}= 136$. The approximating curve is constructed with
constants found from the following two points: $\omega_{1} =2;
\omega_{2}= 30, \alpha=0.5.$ Mean quadratic error: $\sigma =
0.867221.$}
\end{figure}

\begin{figure}[h] 
\centering\epsfig{figure=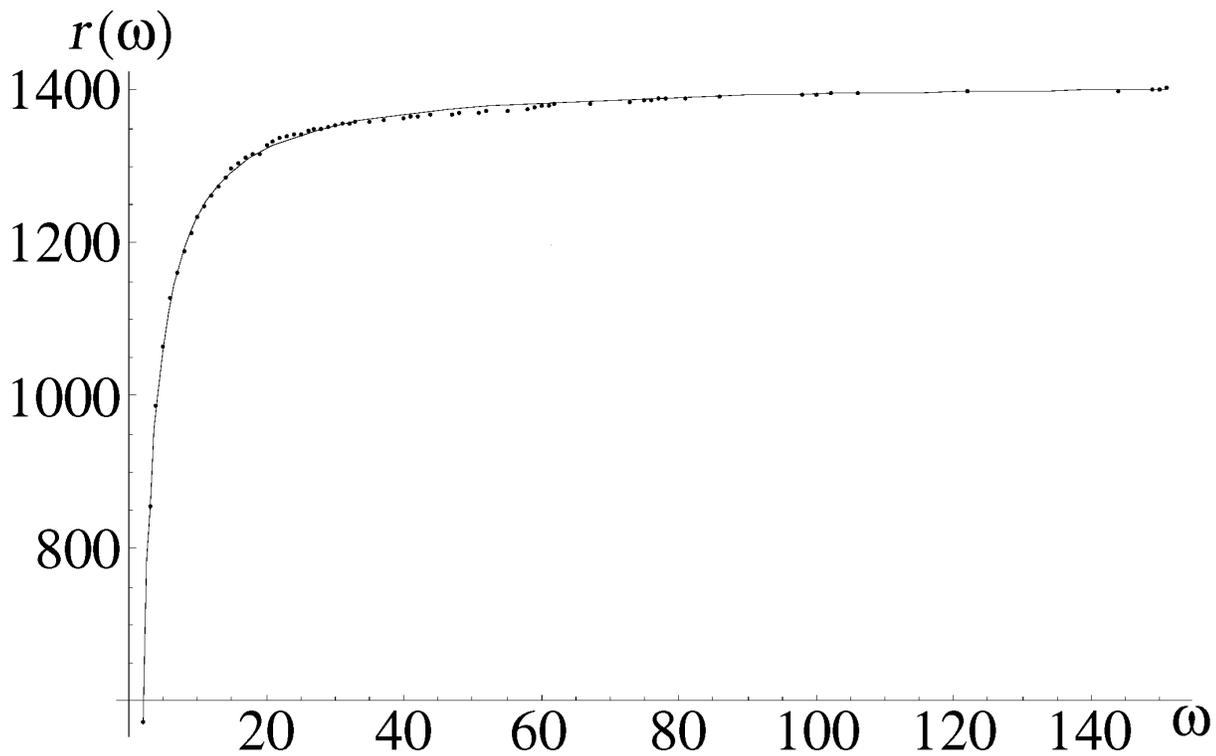,width=\linewidth} \caption{\small
A.N.Ostrovsky "Good Friends Settle Things." [Svoi ljudi -
sochtemsja.] Dependence of the inverted rank on the frequency.
Practical and theoretical dependencies. At the bottom, 2
frequencies are cut off: $\omega_{\min} =1$. At the top, 15
frequencies are cut off: $\omega_{\max}= 151$. The approximating
curve is constructed with constants found from the following two
points: $\omega_{1} =2; \omega_{2}= 30, \alpha=3.4.$ Mean
quadratic error: $\sigma = 0.456089.$}
\end{figure}

\begin{figure}[h] 
\centering\epsfig{figure=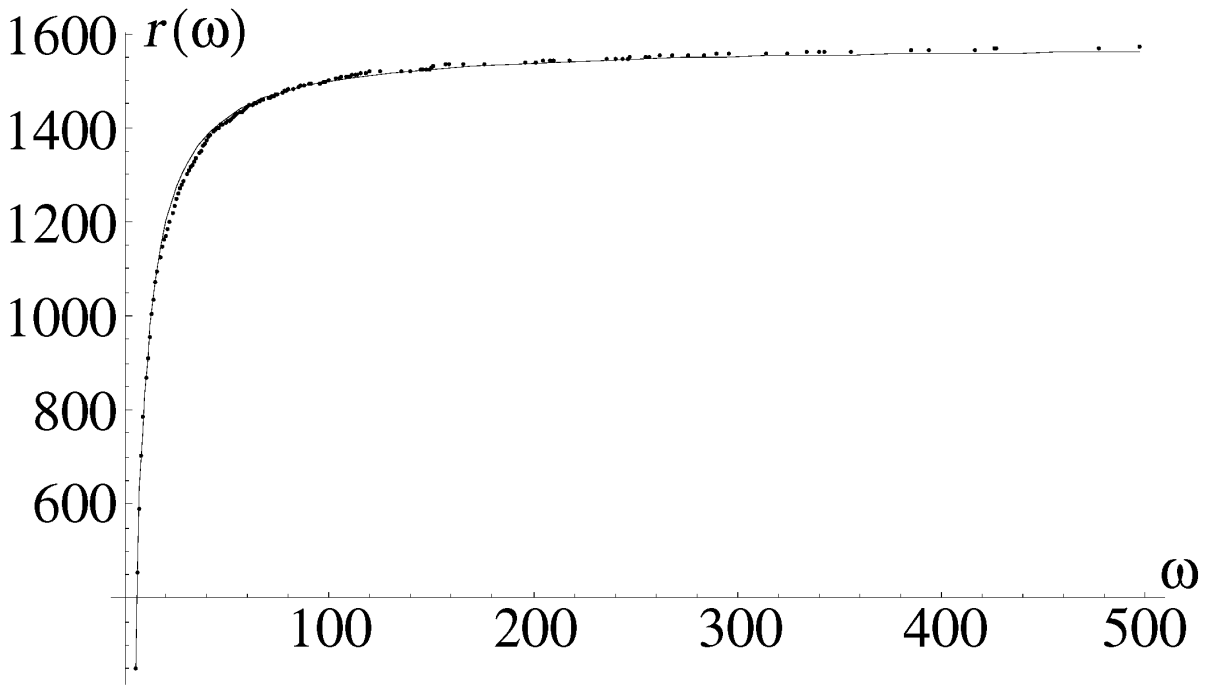,width=\linewidth} \caption{\small
A.N.Ostrovsky. Dictionary for all texts. Dependence of the
inverted rank on the frequency. Practical and theoretical
dependencies. At the bottom, 5 frequencies are cut off:
$\omega_{\min} =1$. At the top, 15 frequencies are cut off:
$\omega_{\max}= 497$. The approximating curve is constructed with
constants found from the following two points: $\omega_{1} =5;
\omega_{2}= 86, \alpha=0.5.$ Mean quadratic error: $\sigma =
1.33802.$}
\end{figure}

\begin{figure}[h] 
\centering\epsfig{figure=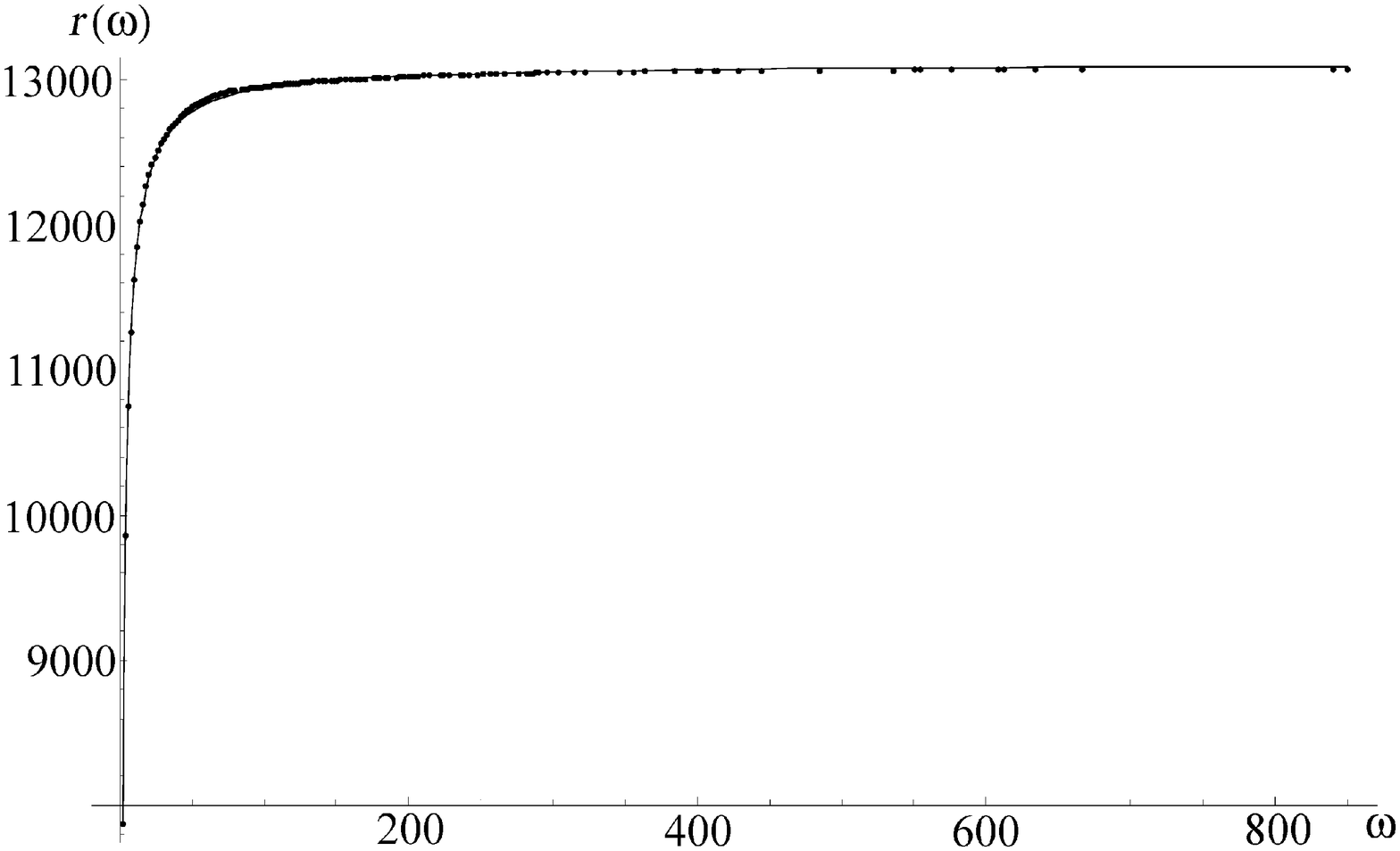,width=\linewidth} \caption{\small
Nicolai Berdyaev. The Fate of Russia. Dependence of the inverted
rank on the frequency. Practical and theoretical dependencies. At
the bottom, 0 frequencies are cut off: $\omega_{\min} =2$. At the
top, 7 frequencies are cut off: $\omega_{\max}= 850$. The
approximating curve is constructed with constants found from the
following two points: $\omega_{1} =List; \omega_{2}= 146,
\alpha=0.35.$ Mean quadratic error: $\sigma = 17.365.$}
\end{figure}

\begin{figure}[h] 
\centering\epsfig{figure=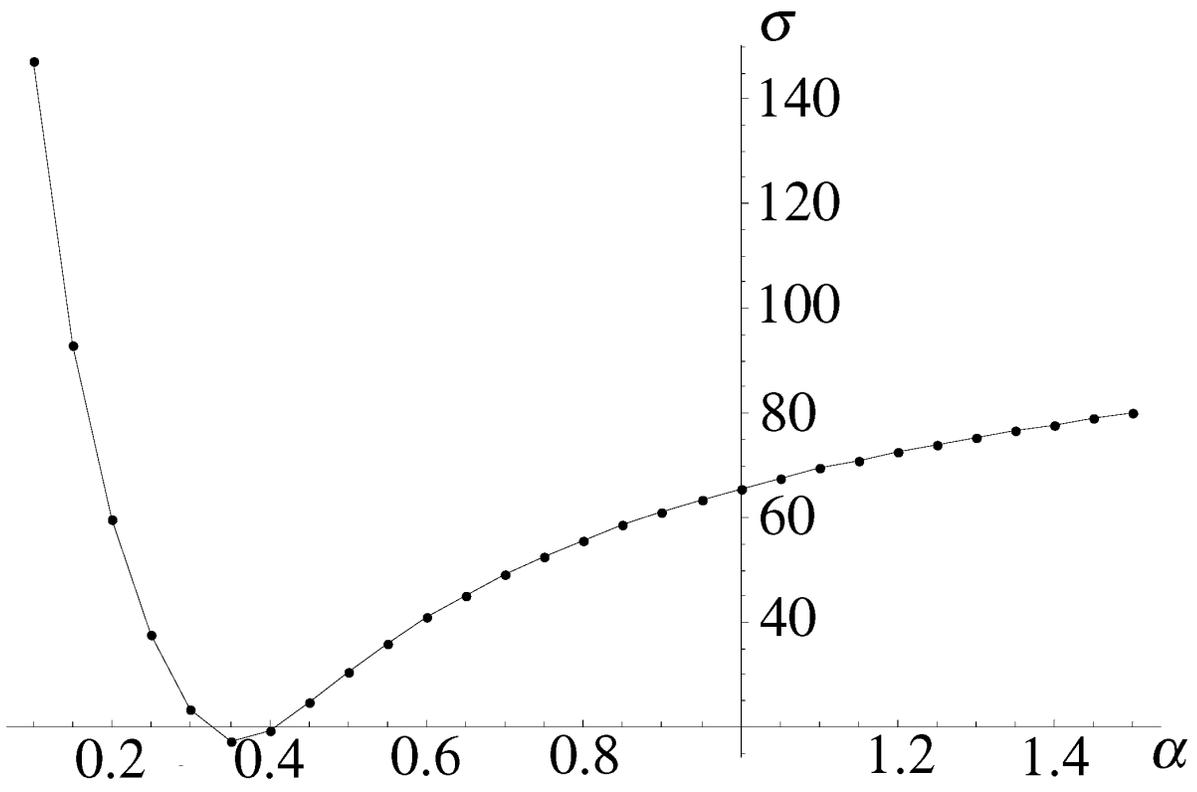,width=\linewidth} \caption{\small
Nicolai Berdyaev. The Fate of Russia. Dependence of the variance
$\sigma$ on the parameter $\alpha$. Range [0.1,1.5] with step
$\delta = 0.05$}.
\end{figure}

\end{document}